\documentclass[11pt]{amsart}
\usepackage{amssymb}
\numberwithin{equation}{section}

\newtheorem{theorema}{Theorem}
\newtheorem{theorem}{Theorem}[section]
\newtheorem{Theorem}[theorem]{Theorem}
\newtheorem{corollary}[theorem]{Corollary}
\newtheorem{proposition}[theorem]{Proposition}
\newtheorem{lemma}[theorem]{Lemma}
\theoremstyle{definition}
\newtheorem{definition}[theorem]{Definition}
\theoremstyle{remark}
\newtheorem{remark}[theorem]{Remark}
\newtheorem{example}[theorem]{Example}

\newcommand\E{\mathcal{E}}

\renewcommand{\O}{\mathcal{O}}

\renewcommand{\L}{\mathcal{L}}

\newcommand{\R}{\mathbb{R}}
\newcommand{\C}{\mathbb{C}}

\newcommand\lie[1]{\mathfrak{#1}}

\newcommand{\fh}{\lie{h}}
\newcommand{\fg}{\lie{g}}
\newcommand{\fm}{\lie{m}}

\newcommand{\fl}{\lie{l}}

\newcommand{\fk}{\lie{k}}

\def    \inv    {^{-1}}

\begin{document}

\title{Contact Reduction}
\author{Christopher Willett}
\address{Department of
Mathematics, University of Illinois, Urbana, IL 61801}
\email{cwillett@math.uiuc.edu}
\date{\today}
\thanks{The author was supported by a National Science 
Foundation graduate Vertical Integration of Research and 
Education  fellowship and the American Institute of Mathematics}

\begin{abstract} In this article I propose a new method for 
reducing a co-oriented contact manifold $M$ equipped with 
an action of a 
Lie group $G$ by contact transformations.
  With a certain regularity and integrality assumption the 
contact quotient $M_\mu$ 
at $\mu \in \fg^*$ is a naturally a co-oriented
contact orbifold which is independent of the contact form 
used to represent the given contact structure.  Removing the regularity 
and integrality
assumptions and replacing them with one concerning the existence of 
a slice, which is satisfied for compact symmetry groups, results in a 
contact stratified space; i.e, a stratified space equipped with a 
line bundle which, when restricted to each stratum, defines a 
co-oriented contact structure.  This extends the previous work of the 
author and E. Lerman (\cite{LW}).
\end{abstract}
\maketitle

\tableofcontents

\section{Introduction} 
In this article I propose a new method for reducing
 a co-oriented contact manifold $M$ under the action of a 
Lie group $G$ by contact transformations.
  With a certain regularity and integrality assumption (see 
Lemma \ref{trans_free} and Theorem \ref{sred})
the contact quotient $M_\mu$ 
at $\mu \in \fg^*$ is a naturally a co-oriented
contact orbifold which is independent of the contact form 
used to represent the given contact structure.  Removing the regularity
and integrality
assumption and replacing it with one concerning the existence of 
a slice, which is satisfied for compact symmetry groups, results in a 
contact stratified space; i.e, a stratified space equipped with a 
line bundle which, when restricted to each stratum, defines a 
co-oriented contact structure.  This extends the previous work of the 
author and E. Lerman (\cite{LW}).

The earliest notion of a contact quotient that I am aware of is in the 
work of V. Guillemin and S. Sternberg, where it appears in the 
guise of co-isotropic 
reduction of symplectic cones (\cite{GS}).  The Guillemin-Sternberg
quotient is based on the reduction of foliations and is akin to the 
Kazhdan-Kostant-Sternberg
quotient in symplectic geometry.  The new contact quotient
 is akin the Marsden-Weinstein-Meyer 
quotient.  The two models differ, but are related:  The Guillemin-Sternberg
quotient fibers over a co-adjoint orbit with typical fiber the new quotient.
The Guillemin-Sternberg quotient is described in a later 
section.  Independently of Guillemin-Sternberg, C. Albert developed 
a model for contact reduction which is valid at all $\mu \in \fg^*$ but
 depends 
upon the contact form used to represent the given 
contact structure (\cite{Al}).

\subsection*{Acknowlegements} The author would like to thank Eugene Lerman 
for 
suggesting this line of investigation and innumerable discussions on the 
material while at the American Institute of Mathematics during 
the fall of 2000
and for reading several preliminary versions of the paper.  The author 
would also like to thank Susan Tolman for making several helpful suggestions 
concerning the presentation of the material.

\subsection*{A Note on Notation}  The notation in this paper will be 
consistent with that used in \cite{LW}.

\section{Group actions on contact manifolds}

Before proceeding to the proof of the main theorems of the paper, we 
first recall the relevant and essential facts about contact 
manifolds and group actions on contact manifolds.

\begin{definition}
A {\bf contact structure} on a manifold $M$ of dimension $2n+1$ is a 
co-dimension one distribution $\xi$ on $M$ which is given locally by 
the kernel 
of a 1-form $\alpha$ with $\alpha \wedge (d \alpha)^n \neq 0$.
If there is a global 1-form $\alpha$ with
$\text{ker }\alpha=\xi$, then $\xi$ is called {\bf co-orientable},
$\alpha$ is said to {\bf represent} $\xi$, and is called a {\bf contact form}.
\end{definition}

\begin{remark} One can define a co-orientable contact structure on a manifold $M$ in the 
following convenient way, which we exploit in a later section.
A co-dimension one distribution $\xi$ is a co-orientable contact 
structure if its annihilator $\xi^\circ$ admits a nowhere vanishing section, 
$\alpha$, and $\xi^\circ -\{0\}$ is a symplectic submanifold of 
$T^*M$.  Equivalently, a line sub-bundle $\zeta \subseteq T^*M$ which is 
symplectic away from the zero section and 
equipped with a nowhere vanishing section defines a co-orientable 
contact structure on $M$.  A choice of a component of $\xi^\circ-\{0\}$ 
is a {\bf co-orientation} for $\xi$.
\end{remark}

If $f$ is a non-vanishing function on $M$, then $f \alpha$ is 
another contact form on $M$ which represents $\xi$.  The 
{\bf conformal class} of $\alpha$ is 
$$\{f \alpha|f \text{ is a positive function }\}$$
It should be 
stressed that the important data on a contact manifold is the contact 
structure, not the contact form used to represent it.
Observe that for each $m \in M$, 
$$(\text{ker }\alpha_m,d \alpha_{m|\text{ker }\alpha_m})$$ is a 
symplectic vector space.  On a co-oriented contact manifold $(M,\alpha)$
 there is a distinguished vector field $Y$, 
called the {\bf Reeb} vector field, which satisifies
$$\iota(Y)\alpha=1 \hskip.25in \iota(Y) d \alpha=0$$
$Y$ is unique with respect to these two properties.
This allows us to split the tangent bundle of $M$ as 
$$TM=\text{ker }\alpha \oplus \R Y$$ where $\R Y$ is the line bundle over $M$ 
spanned by $Y$.
Since $\L_Y \alpha =d \iota(Y)d \alpha+d \iota(Y) \alpha=0$, the 
 flow of the Reeb vector field preserves $\alpha$.

For the purposes of this article, the most important examples of contact 
manifolds come via interaction with symplectic manifolds. 

\begin{definition} A hypersurface $\Sigma$ of a symplectic
manifold, $(N,\omega)$ is said to be of {\bf contact type}
if there is a vector field $X$, defined near $\Sigma$,
 which satisfies $\L_X \omega=\omega$.  $X$ is called a {\bf Liouville} vector
field.
\end{definition}
If $\Sigma$ is a hypersurface of contact type in a symplectic manifold $(N,\omega)$
with Liouville vector field $X$, then 
$(\Sigma,\iota(X)\omega_{|\Sigma})$ is a contact manifold.  

\begin{example}  
If $(M,\xi)$ is a contact manifold, then away from the 
zero section the annihilator of $\xi$ in 
$T^*M$
is a symplectic submanifold of $T^*M$.
A connected component of $\xi^\circ -\{0\}$ 
is called the {\bf symplectization} of $(M,\xi)$.
If $(M,\alpha)$ is a co-oriented contact manifold, then the symplectization can
be identified with 
$(M \times \R,d(e^t \alpha))$, where $t$ is the $\R$ coordinate.  Note that every 
co-oriented 
contact manifold is a hypersurface of contact type in its symplectization.
\end{example}

\begin{definition}  Suppose a Lie group $G$ acts on a contact 
manifold $(M,\alpha)$ and preserves $\alpha$.
The {\bf contact moment map associated} to the action
 is denoted by $\Phi_\alpha:M \to \fg^*$ and defined by 
$$\langle \Phi_\alpha(x),A \rangle=\alpha_x(A_M(x))$$ for $x \in M$ and 
all $A \in \fg$. The contact moment map is equivariant, where $G$ acts on 
$\fg^*$ through the co-adjoint respresentation (\cite{Geiges}).
\end{definition}

\begin{remark}
The moment map depends upon the choice of an invariant 
contact form to represent the contact structure:
If $f$ is a non-vanishing invariant function on $M$, then 
$f \alpha$ is an invariant contact form on $M$ 
and $\Phi_{f \alpha}=f \Phi_\alpha$. 
\end{remark}
 
\begin{example}   Suppose $(N,\omega)$ is a symplectic 
manifold equipped with a Hamiltonian $G$-action and let $\Psi:N \to \fg^*$
be a corresponding 
equivariant symplectic 
moment map.  Let $\Sigma \subset N$ be an 
invariant hypersurface which has an invariant Liouville vector field and 
$\alpha$ the induced invariant contact form on $\Sigma$. Then the 
restriction of the $G$-action to $\Sigma$ preserves $\alpha$
and the associated
contact moment map $\Phi_\alpha:\Sigma \to \fg^*$ is the restriction of $\Psi$ to 
$\Sigma$.
\end{example}

\begin{example}  Suppose a Lie group $G$ acts on a contact manifold $(M,\alpha)$
preserving $\alpha$ and let $\Phi_\alpha:M \to \fg^*$ be the associated contact 
moment map.  The extension of the action to the symplectization $(M \times \R, 
d(e^t \alpha))$  by $g \cdot(m,t)=(g \cdot m,t)$ is Hamiltonian and a corresponding
equivariant symplectic moment map $\Psi:M \times \R \to \fg^*$ is given by 
$\Psi(m,t)=e^t \Phi_\alpha(m)$.
\end{example}

\begin{lemma} \label{inv_subspace}
Suppose a Lie group $G$ acts on a contact manifold $(M,\alpha)$, 
preserving $\alpha$, and let $\Phi_\alpha:M \to \fg^*$ be the 
associated moment map.  Let $\Theta \subset \fg^*$ be a subset which is 
invariant under dilations by $\R^+$.
\begin{itemize}
\item Let $\Psi:M \times \R \to \fg^*$ be the symplectic moment map for the
induced action on $M \times \R$. Then $\Psi\inv(\Theta)=\Phi_\alpha\inv(\Theta) 
\times \R$.
\item  If $f$ is a positive invariant function on $M$, then 
$\Phi_{f \alpha}\inv(\Theta)=\Phi_\alpha\inv(\Theta)$.
\end{itemize}
\end{lemma}
\begin{proof}
Because  $\Psi(m,t)=e^t \Phi_\alpha(m)$ and $\Theta$ is $\R^+$ 
invariant, $(m,t) \in \Psi\inv(\Theta)$ if and only if 
$\Phi_\alpha(m) \in \Theta$.  This establishes the first point. 
The second one follows in the same manner, using the fact that 
$\Phi_{f \alpha}=f \Phi_\alpha$.
\end{proof}

Recall that an action of a Lie group $G$ on a manifold $M$ is 
called {\bf proper} if the map $G \times M \to M \times M$ given by 
$(g,m) \mapsto (g \cdot m, m)$ is a proper map.  Actions of compact 
groups are proper.  Proper actions have compact isotropy groups.
If a Lie group $G$ acts properly
 on a (paracompact) co-oriented contact manifold $M$ and preserves 
the contact structure and a co-orientation,
 then we can find an invariant contact form, $\alpha$, on 
$M$ which 
represents 
the contact structure (\cite{IJL}). In this case, the Reeb vector field is, 
by uniqueness, $G$-invariant as well.

A {\bf slice} for an action of a Lie group $G$ on a manifold $M$ at a
point $x$ is a $G_x$ invariant submanifold $S$ such that $G\cdot
S$ is an open subset of $M$ and such that the map $G\times S
\to G\cdot S$, $(g, s) \mapsto g\cdot s$  descends to 
a diffeomorphism $G\times _{G_x} S \to G\cdot S$, $[g, s]
\mapsto g\cdot s$.  A
theorem of Palais  asserts that for smooth proper actions
slices exist at every point (\cite{Palais}).

\begin{proposition}\label{contact_cross_section}  
Suppose a Lie group $G$ acts properly
on a contact manifold $(M,\alpha)$, preserving $\alpha$, 
and let $\Phi_\alpha:M \to \fg^*$ 
be the associated contact moment map. 
If there exists a slice $S$ through $x \in \fg^*$ for the 
co-adjoint action which is 
invariant under dilations by $\R^+$, then
 $R:=\Phi_\alpha\inv(S)$ is a $G_x$ invariant contact submanifold of 
$(M,\alpha)$ and 
the contact moment map for the $G_x$ action on $R$ is given by the 
restriction of the 
$\Phi_\alpha$ to $R$ followed by the natural projection of $\fg^*$ 
onto $\fg_x^*$.
$R$ is called a {\bf contact cross section}.
\end{proposition}

\begin{proof} 

The proof of Proposition \ref{contact_cross_section} 
is derived from 
the corresponding proof of the symplectic cross section theorem of 
Guillemin and Sternberg (\cite{convex}).  Indeed, extend the $G-$action 
to the symplectization of $M$ and let $\Psi$ be the corresponding symplectic 
moment map.  The symplectic cross section theorem gives $\Psi\inv(R)$ as a 
symplectic submanifold of the symplectization. By hypothesis, 
$S$ is invariant under dilations by $\R^+$ and hence Lemma \ref{inv_subspace} 
implies that $\Psi\inv(R) = \Phi_\alpha\inv(R) 
\times \R$.
It follows that
 the contact cross section is a hypersurface of contact type in the 
symplectic cross section, which gives the result.
\end{proof}

\begin{remark} Although Proposition \ref{contact_cross_section} can be 
proved directly, the direct proof does not diverge significantly from the 
proof of the symplectic cross section theorem of Guillemin and Sternberg 
found in \cite{convex}.
\end{remark}

\section{The Reduction Theorems}
In this section we prove the first two of our three main theorems.

\begin{proposition}\label{basics}  Suppose a Lie group $G$ acts 
properly on a contact 
manifold $(M,\alpha)$, preserving $\alpha$ and let 
$\Phi_\alpha:M \to \fg^*$ be the associated moment map. 
\begin{enumerate}

\item The Reeb flow preserves the level sets of $\Phi_\alpha$.

\item For all $x \in M, \vec v \in T_xM,$ and $A \in \fg$, 
$$\langle d (\Phi_\alpha)_x(\vec v),A \rangle=d \alpha_m(\vec v,A_M(x))$$

\item If $\Phi_\alpha(x)=0$, then $T_x(G \cdot x)$ is an isotropic subspace of the 
symplectic vector space $(\text{ker }\alpha_x,d \alpha_x)$.

\item $\text{Im}(d (\Phi_\alpha)_x)^\circ=\{A \in \fg|\iota(A_M(x))d \alpha_x=0\}=
\{A \in \fg|A_M(x) \in \text{ker }d \alpha_x\}$
\end{enumerate}
\end{proposition}

\begin{proof}  Denote the Reeb vector field by $Y$ and its 
flow by $\rho_t$.  
As noted earlier, by uniqueness $Y$ is invariant and 
therefore $\rho_t$ is 
equivariant.  Therefore, for any $A \in \fg$ and $m \in M$, 
it follows that $A_M(\rho_t(m))=d \rho_t(A_M(m))$.  Thus 
$$\aligned 
\langle \Phi_\alpha(\rho_t(z)),A \rangle &=\alpha_{\rho_t(z)}(A_M(\rho_t(z))) \cr
&=\rho_t^*(\alpha)_z(A_M(z)) \cr
&=\alpha_z(A_M(z)) \cr
&=\langle \Phi_\alpha(z),A \rangle 
\endaligned 
$$
This establishes item (1).

Since $\L_{A_M} \alpha=0$, 
Cartan's formula gives $d \iota(A_M) \alpha=-d \alpha(A_M,-)$.
Therefore, 
$$\langle d (\Phi_\alpha)_x(\vec v),A \rangle=d \langle \Phi_\alpha,
A \rangle_x(\vec v)
=d \alpha_x(\vec v,A_M(x))$$ which establishes item (2).

The proof of item (3) is an application of item (2).
For any $A,B \in \fg$, the second item implies that for any $x \in \Phi_\alpha\inv(0)$, 
$$\aligned d \alpha_x(A_M(x),B_M(x))&=\langle d(\Phi_\alpha)_x(A_M(x)),B \rangle \cr
&=\langle A_{\fg^*}(\Phi_\alpha(x)),B \rangle \hskip.25in \text{by equivariance}\cr
&=\langle A_{\fg^*}(0),B \rangle \cr
&=0
\endaligned$$
Thus, $A \in T_x(G \cdot x)^{d \alpha_x}$, whence $T_x(G \cdot x) \subseteq T_x(G \cdot x)^{d \alpha_x}$.

The final point follows immediately from the second point.
Note that $\text{ker }d \alpha_x=\R Y(x)$.
\end{proof}

\begin{remark} Note that the last item of Proposition \ref{basics} 
implies that $x \in M$ is a regular point for $\Phi$ if and only if 
$$\text{dim }\{A \in \fg|A_M(x) \in \R Y(x)\}=0$$
That is, $x \in M$ may be critical yet 
have a discrete stabilizer. 
 This contrasts sharply with the symplectic category, where 
a point is regular if and only if its stabilizer is discrete.
\end{remark}

\begin{lemma} \label{lie_ideal}  Let $G$ be a Lie group and choose
$\mu \in \fg^*$.  Then $\fk_\mu=\text{ker }\mu \cap \fg^*=
\text{ker}(\mu_{|\fg_\mu})$
is a Lie ideal in $\fg_\mu$.  Hence, there is a unique, connected 
normal Lie subgroup of $G_\mu$ with Lie algebra $\fk_\mu$.
\end{lemma}

\begin{proof}
Choose $A \in \fg_\mu$.  Then $0=A_{\fg^*}(\mu)=\text{ad}^\dagger(A)\mu$,
where $\text{ad}^\dagger:\fg \to \text{End}(\fg^*)$ is the differential 
of the co-adjoint representation.
Thus, for any $B \in \text{ker }(\mu_{|\fg_\mu})$, we have that 
$$0=\langle \text{ad}^\dagger(A)\mu,B \rangle=\langle \mu,[A,B] \rangle$$
Therefore, $\fk_\mu$ is a Lie ideal in $\fg_\mu$.
\end{proof}

\begin{definition}  Suppose a Lie group $G$ acts properly on a 
contact manifold, $(M,\alpha)$, preserving $\alpha$, and 
let $\Phi_\alpha:M \to \fg^*$ be the associated moment map.
Choose $\mu \in \fg^*$.  We define the {\bf kernel group} of 
$\mu$ to be the unique connected 
Lie subgroup of $G_\mu$ with Lie algebra 
$\fk_\mu:=\text{ker }\mu_{|\fg_\mu}$. The kernel group of $\mu$ is denoted by 
$K_\mu$.  We define
the {\bf contact quotient} (or {\bf contact 
reduction}) of $M$ by $G$ at $\mu$ to be 
$$M_\mu:=\Phi_\alpha\inv(\R^+\mu)/K_\mu$$
\end{definition}

\begin{remark}
If $f$ is a positive invariant function on $M$, then $f \alpha$ is another 
invariant contact 
form on $M$ which is in the same conformal class as $\alpha$.  
Since $\R^+ \mu$ is invariant under dilations by $\R^+$,
Lemma \ref{inv_subspace} implies 
 that $\Phi_{f \alpha}\inv(\R^+ \mu)=\Phi_{\alpha}\inv(\R^+\mu)$ and hence the reduced space is topologically independent of the choice of contact forms in the same 
conformal class.
\end{remark}

\begin{lemma}\label{trans_free}
Suppose a Lie group $G$ acts on a contact manifold $(M,\alpha)$, preserving 
$\alpha$, and let $\Phi_\alpha:M \to \fg^*$ be the associated contact 
moment map.  Choose $\mu \in \fg^*$ and let $K_\mu$ be the connected 
Lie subgroup of $G_\mu$ with Lie algebra $\fk_\mu=\text{ker }\mu_{|\fg_\mu}$.
Then $\Phi_\alpha$ is transverse to $\R^+ \mu$ if and only if $K_\mu$ 
acts locally freely on $\Phi_\alpha\inv(\R^+ \mu)$.

\end{lemma}
\begin{proof}
 If $\Phi_\alpha$ is transverse to $\R^+ \mu$, then
$$Z:=\Phi_\alpha\inv(\R^+ \mu)$$
is a submanifold of $M$.
  Choose $z \in Z$.  The transversality 
condition, $$\text{Im}(d (\Phi_\alpha)_z) + \R \mu = \fg^*$$ is equivalent 
to the condition $$\text{Im}(d (\Phi_\alpha)_z)^\circ \cap \text{ker }\mu=0$$
Let $H$ be the isotropy subgroup of $z$ in $K_\mu$ and 
$\fh$ its Lie algebra.
 Choose $A \in \fh$ and denote the Reeb vector field of $(M,\alpha)$ by $Y$.
 Then $\langle \mu,A \rangle =0$ since 
$\fh \subseteq \fk_\mu$.  Because $H$ fixes $z$, $A_M(z)=0$.  Thus,
$$\aligned 
A & \in \{B \in \fg | B_M(z) \in \R Y(z) \} \cap \text{ker }\mu \cr
&=\text{Im}(d (\Phi_\alpha)_z)^\circ \cap \text{ker }\mu \hskip.25in 
\text{by Proposition }
\ref{basics} \cr
&=0
\endaligned$$
Therefore, $\fh=0$ and $K_\mu$ acts locally freely on $Z$.

To establish the second statement, suppose that 
$\Phi_\alpha(z)=s \mu$ for some $s \in \R^+$ and let
$A \in \text{ker }\mu \cap (\text{Im}(d(\Phi_\alpha)_z))^\circ$.  Then,
by Proposition, \ref{basics}
$A_M(z)=t Y(z)$ for some $t \in \R$, where $Y$ is the Reeb vector
field.
Therefore,
$$0=\langle s \mu,A \rangle = \langle \Phi_\alpha(z),A \rangle=\alpha_z(A_M(z))
=\alpha_z(t Y(z))=t$$
Hence, $A_M(z)=0$ and $A \in \fg_z$. Since $A \in \text{ker }\mu$ and 
$A \in \fg_z \subseteq \fg_\mu$, it follows that $A \in \fk_\mu$.
Since $K_\mu$ acts locally freely on $\Phi_\alpha\inv(\R^+ \mu)$, 
$A=0$ and hence $\Phi_\alpha$ is transverse to $\R^+ \mu$.

\end{proof}

While the contact quotient is defined at any element of $\fg^*$, it 
is not necessarily a contact manifold.  There are two problems.
The first is that there is no guarantee that the kernel group of 
$\mu$ will act properly on $\Phi\inv(\R^+\mu)$, in which case 
the resulting quotient may not be Hausdorff.
If $G$ is compact 
and $\mu$ is integral (i.e, $K_\mu=\text{ker }\chi_\mu$, where 
$\chi_\mu:G_\mu \to S^1$ is a group map satisfying 
$d \chi_u=\mu_{|\fg_\mu}$), 
then $K_\mu$ is actually compact and no hypothesis is needed.
This integrality condition is required in the Guillemin-Sternberg 
procedure, described later (\cite{GS}).
The second problem is that the kernel and isotropy groups of 
$\mu$ may coincide.  If $\mu$ is non-zero, then the resulting 
quotient may fail to be contact, as the below example 
shows.  Hence, we assume that $\text{ker }\mu
+\fg_\mu=\fg$.  If $G$ is compact, then the existence of an invariant 
metric on $G$ implies this condition.

\begin{example} \label{sl2_example} Let $G=\text{SL}(2,\R)$ and let 
$\alpha$ be the natural contact form on $M=T^*(G) \times \R \cong
G \times \fg^* \times \R$.  For the lift of the natural $G$ action the 
associated contact moment map, $\Phi_\alpha:M \to \fg^*$, is given by 
$\Phi_\alpha(g,\vec v,t)=\text{Ad}^\dagger(g)\cdot \vec v$.  Set 
$$\mu=
\left( \begin{array}{ccc}
0 & 1 \\
0 & 0
\end{array}\right)
$$
Then 
$$\fk_\mu=\fg_\mu=\{\left( \begin{array}{ccc}
0 & t \\
0 & 0
\end{array}\right)|t \in \R\}
$$
Therefore, $$K_\mu=\{\left( \begin{array}{ccc}
1 & t \\
0 & 1
\end{array}\right)|t \in \R\}
$$
is closed in $G$ and acts properly on $M$, but
 $\Phi_\alpha\inv(\R^+\mu)/K_\mu$ is four dimensional and hence isn't contact.
\end{example}

\begin{lemma} \label{kernel_splitting}  Let $V$ be a vector space and 
$\omega:V \times V \to \R$ be an antisymmetric bilinear map.  Suppose there is a 
decomposition $V = X \oplus W$ which is perpendicular with respect to 
$\omega$; i.e, $\omega(x,w)=0$ for all $x \in X$ and $w \in W$.  If 
$\text{ker }\omega \subseteq \text{ker }\omega_{|X}$, then 
$\text{ker }\omega = \text{ker }\omega_{|X}$.
\end{lemma}

\begin{proof}  Choose $x_0 \in \text{ker }\omega_{|X}$ and $v \in V$. 
Write $v = x+w$ where $x \in X$ and $w \in W$.  Then 
$\omega(x_0,v)=\omega(x_0,x)+\omega(x_0,w)=0$.  Hence $x_0 \in \text{
ker }\omega$.
\end{proof}

\begin{lemma} \label{isotrop_ker}  Let $(V,\omega)$ be a symplectic 
vector space and suppose that $W$ is an isotropic subspace.  Then 
$$\text{ker }\omega_{|W^\omega}=W$$
where $W^\omega$ is the symplectic perpendicular of $W$ with respect to 
$\omega$.
\end{lemma}

\begin{proof}  Fix $w \in W$.  Since $W$ is isotropic, $w \in W^\omega$.
For any $x \in W^\omega$, $\omega(w,x)=0$ by definition.  Hence 
$W \subseteq \text{ker }\omega_{|W^\omega}$.  Conversely, if 
$x \in \text{ker }\omega_{W^\omega}$, then by unravelling the 
various definitions, one sees that $x \in (W^\omega)^\omega=W$.
\end{proof}

\begin{theorema}\label{sred} 
Suppose a Lie group $G$ acts 
 on 
a contact manifold, $(M,\alpha)$, preserving $\alpha$, and let 
$\Phi_\alpha:M \to \fg^*$ be the associated contact moment map.  Choose 
$\mu \in \fg^*$ and let $K_\mu$ be the connected 
Lie subgroup of $G_\mu$ with Lie algebra $\fk_\mu=\text{ker }\mu_{|\fg_\mu}$.
If 
\begin{itemize}
\item $K_\mu$ acts properly on $\Phi_\alpha\inv(\R^+ \mu)$
\item $\Phi$ is transverse to $\R^+ \mu$
\item $\text{ker }\mu + \fg_\mu=\fg$
\end{itemize}
then 
 the quotient $$M_\mu=\Phi_\alpha\inv(\R^+\mu)/K_\mu$$ is naturally 
a contact orbifold; i.e, $$\text{ker }\alpha \cap T(\Phi_\alpha\inv(\R^+ \mu))$$
descends to a contact structure on the quotient.
\end{theorema}

\begin{remark}  In the special case of $\mu=0$
C. Albert, H. Geiges, and F. Loose, independently of Guillemin-Sternberg, 
established the above theorem in various papers (\cite{Al, Geiges, 
Loose}).
\end{remark}

\begin{proof}
Since $\Phi_\alpha$ is transverse to $\R^+ \mu$, 
$Z:=\Phi_\alpha\inv(\R^+ \mu)$ is a submanifold of $M$ and
Lemma \ref{trans_free} implies that $K_\mu$ acts locally freely $Z$.
Hence, $M_\mu$ is an orbifold.  Fix $z \in Z$.
For any $A \in \fk_\mu$, we have that
$$\aligned 
(\iota(A_M)\alpha)(z)&=\langle \Phi_\alpha(z),A \rangle \cr
&=\langle s \mu,A \rangle \hskip.25in \text{for some }s \in \R^+ \cr
&=0
\endaligned $$
  Hence, $\alpha$ descends to $\alpha_\mu$ on
$M_\mu$, which is contact if and only if 
$$\text{ker}(d \alpha_{z|T_zZ \cap \text{ker }\alpha_z})=T_z(K_\mu \cdot z)$$
Let $\Psi_\alpha:M \to \fk_\mu^*$ be the contact moment map associated 
to the action of $K_\mu$ on $M$.  Note that $\Psi_\alpha=i^T \circ \Phi_\alpha$
where $i^T:\fg^* \to \fk_\mu^*$ is the natural projection.
Observe that $\Psi_\alpha(z)=0$ and hence that
$T_z(K_\mu \cdot z)$ is an isotropic subspace of the symplectic 
vector space $(\text{ker }\alpha_z,d \alpha_{z|\text{ker }\alpha_z})$
 by Proposition \ref{basics}.
Let $T_z(K_\mu \cdot z)^{d \alpha_z}$ denote the symplectic perpendicular
of $T_z(K_\mu \cdot z)$ in $(\text{ker }\alpha_z,d \alpha_{z|\text{ker }\alpha_z})$.
A vector $\vec v$ is in  $T_z(K_\mu \cdot z)^{d \alpha_z}$ if and only if for all 
$A \in \fk_\mu$
$$\aligned
0&=d \alpha_z(\vec v, A_M(z)) \cr
&=\langle d (\Phi_\alpha)_z(\vec v),A \rangle \hskip.25in 
\text{by Proposition } \ref{basics} 
\endaligned$$
That is, $\vec v \in T_z(K_\mu \cdot z)^{d \alpha_z}$ if and only if 
$d (\Phi_\alpha)_z(\vec v)_{|\fk_\mu} \equiv 0$. Hence, 
$$T_z(K_\mu \cdot z)^{d \alpha_z}=T_z U \cap \text{ker }\alpha_z$$ where 
$U=\Phi_\alpha\inv(\fk_\mu^\circ)$.  Note that $U$ is a submanifold of 
$M$ by the transversality condition. Indeed, $U=\Phi_\alpha\inv(\fk_\mu^\circ)
=\Psi_\alpha\inv(0)$.  Since $\Phi_\alpha$ is transverse to 
$\R^+ \mu$, $K_\mu$ acts locally freely on $Z$ by Lemma \ref{trans_free}
and thus on a neighborhood of $Z$ and hence on $U$ (or, at least on a 
neighborhood of $Z$ in $U$).
  Therefore, also by Lemma \ref{trans_free},
$U$ is a submanifold of $M$.  Lemma \ref{isotrop_ker} implies that
$T_z(K_\mu \cdot z) =\text{ker}(d \alpha_{z|T_zU \cap \text{ker }\alpha_z})$.

We first show that 
$T_z(K_\mu \cdot z) \subseteq \text{ker}(d \alpha_{z|T_zZ \cap \text{ker }\alpha_z})$.
 If $A_M(z) \in T_z(K_\mu \cdot z)$ and $\vec v \in T_zZ \cap \text{ker }
\alpha_z$, then 
$$\aligned d \alpha_z(A_M(z),\vec v) &= -\langle d (\Phi_\alpha)_z(\vec v),
A \rangle \cr
&=s \langle \mu,A \rangle \cr
&=0
\endaligned$$ Therefore, $$T_z(K_\mu \cdot z) 
\subseteq \text{ker}(d \alpha_{z|T_zZ \cap \text{ker }\alpha_z})$$
Note that this implies that 
$\text{ker}(d \alpha_{z|T_zU \cap \text{ker }\alpha_z}) \subseteq 
\text{ker}(d \alpha_{z|T_zZ \cap \text{ker }\alpha_z})$.

The reverse inclusion is slightly more delicate.  Since $\text{ker }\mu + \fg_\mu=\fg$, 
we can choose a splitting $\fg=\fg_\mu \oplus \fm$ where $\mu_{|\fm} \equiv 0$.
Let $\fm_M(z)=\{A_M(z)|A \in \fm \}$.  The proof is completed by showing that 
$\fm_M(z)$ and $T_zZ \cap \text{ker }\alpha_z$ are complementary subspaces of 
$T_zU \cap \text{ker }\alpha_z$ which are 
prependicular with respect to $d \alpha_{z|T_zU \cap \text{ker }\alpha_z}$.
Lemma \ref{kernel_splitting} implies the reverse inclusion.

 We first show that $\fm_M(z) \subset T_zU \cap \text{ker }
\alpha_z$.
Choose any $A \in \fm$ and let $B \in \fk_\mu$.  Then for some $t \in \R^+$,
$$\aligned 
\langle d (\Phi_\alpha)(A_M(z)),B \rangle&=\langle A_{\fg^*}(t \mu),B \rangle 
\hskip.25in \text{ by equivariance}\cr
&=t \langle \mu,[A,B] \rangle \cr
&=-t \langle B_{\fg^*}(\mu),A \rangle \cr
&=0
\endaligned$$
since $B \in \fk_\mu \subset \fg_\mu$.
  Hence, $\fm_M(z) \subset T_zU$.
Additionally,
$$\aligned
\alpha_z(A_M(z))&=\langle \Phi_\alpha(z),A \rangle\cr
&=t \langle \mu,A \rangle \hskip.25in \text{for some }t \in \R^+\cr
&=0 \hskip.5in \text{since }\mu_{|\fm}\equiv 0
\endaligned$$
Thus, $\fm_M(z) \subset T_zU \cap \text{ker }\alpha_z$.

By equivariance, $d(\Phi_\alpha)_z(A_M(z))=A_\fg^*(t \mu)$.  Hence, 
$d(\Phi_\alpha)_z$ maps $\fm_M(z)$ to $T_{t \mu}(G \cdot t \mu)$.  In 
fact, one can show that this map is an isomorphism.
 By definition, $d(\Phi_\alpha)_z(T_zZ)=\R \mu$.
The assumption $\text{ker }\mu+\fg_\mu=\fg$ is equivalent to 
$$\aligned 0&=(\text{ker }\mu)^\circ \cap (\fg_\mu)^\circ \cr 
&=\R \mu \cap T_{t \mu}(G \cdot t \mu)
\endaligned$$
Hence, $\fm_M(z) \cap T_zZ=\{0\}$.  Finally, if $A_M(z)=0$, then 
$A \in \fg_z \subseteq \fm_\mu$.  Hence, $\fm_M(z)$ is a subspace 
of dimension $\text{dim }\fm=\text{dim }\fg-\text{dim }\fg_\mu$.
A dimension count implies that $\fm_M(z)$ and $T_zZ \cap \text{ker }
\alpha_z$ are complementary subspaces of $T_zU \cap \text{ker }\alpha_z$, 
which completes the proof.
\end{proof}

The second main theorem removes the transversality condition above 
but replaces it with a requirement that a convex, $\R^+$
invariant slice exists for $\mu \in \fg^*$.  This condition is 
satisfied for compact $G$.  It is 
suspected, but not known to the author, that a stratification 
theorem for proper group actions ought to hold.  In the 
particular case of $\mu=0$, the author and E. Lerman showed that 
$M_0$ is topologically a stratified space:

\begin{Theorem}   \label{tscsq}
Let $M$ be a manifold with a co-oriented contact structure $\xi$.
Suppose a Lie group $G$ acts properly on $M$ preserving $\xi$ and a 
co-orientation for $\xi$.  Choose
a $G$-invariant contact form $\alpha$ with $\ker \alpha = \xi$ and let
$\Phi_\alpha:M\to \fg^*$ be the corresponding moment map.

Then  for every subgroup $H$ of $G$, each connected component of
the topological space 
$$
\left( M_{(H)} \cap \Phi_\alpha\inv (0) \right)/G
$$ 
is a manifold and the partition of the contact quotient 
$$
M_0 \equiv M/\!/G: = \Phi_\alpha\inv (0)/G 
$$ 
into these manifolds is a stratification.  The symbol
$M_{(H)}$ stands for the set of points in $M$ with the isotropy groups
conjugate to $H$.
\end{Theorem} 

\begin{proof}  See \cite{LW} for a complete proof.
\end{proof}

\begin{proposition}\label{slice_reduction} Suppose a 
 Lie group $G$ acts  on a contact manifold 
$(M,\alpha)$, preserving $\alpha$, and let $\Phi_\alpha:M \to \fg^*$ be the 
associated contact 
moment map. Choose $\mu \in \fg^*$.
Suppose 
\begin{itemize}
\item The kernel group $K_\mu$ of $\mu$ acts properly on $R:=\Phi_\alpha\inv(S)$.
\item There is a convex slice $S$ for $\mu$ which is invariant under 
dilations by $\R^+$
\end{itemize}
 Denote the natural projection,  $\fg^* \to \fg_\mu^*$ by $i^T$ and let
$\mu'=i^T(\mu)$.  Then $M_\mu=R_{\mu'}$, where $R_{\mu'}$ is the 
reduction of $R$ by $G_\mu$ at $\mu'$. 
\end{proposition} 

\begin{proof}  By Proposition \ref{contact_cross_section} 
the restriction of $\alpha$ to $R$ is a contact form and $R$ is 
$G_\mu$ invariant.
The contact moment map for the $G_\mu$ action on 
$R$ is given by $\Psi_\alpha=i^T \circ \Phi_{\alpha|R}$.
Note that $$\Psi_\alpha\inv(\R^+ \mu')=\Phi_\alpha\inv(\R^+ \mu + 
\fg_\mu^\circ)$$  The proof will be completed by showing that 
$(\R^+\mu + \fg_\mu^\circ) \cap S=\R^+ \mu$.  
From this it follows that 
$\Phi_\alpha\inv(\R^+ \mu)=\Psi_\alpha\inv(\R^+ \mu')$.  Since 
$K_\mu=K_{\mu'}$, the conclusion follows.

Since $S$ is $\R^+$-invariant and contains $\mu$, 
$\R^+\mu \subseteq (\R^+\mu + \fg_\mu^\circ) \cap S$.
Recall that $\fg_\mu^\circ =T_\mu(G \cdot \mu)$.
If $t \mu +\eta \in (\R^+ \mu + \fg_\mu^\circ) \cap S$, 
then $\gamma(s)=\mu+\frac{s}{t}\eta, s \in [0,1]$ is tangent to
$G \cdot \mu$ at $\mu$.  On the other hand, since $\gamma(0), \gamma(1) \in S$ 
and $S$ is convex, it follows that $\gamma$ is also tangent to 
$S$ at $\mu$.  Hence $\gamma \in T_\mu (G \cdot \mu) \cap T_\mu S$
However, because $S$ is a slice, $T_\mu(G \cdot \mu) \cap T_\mu S=0$.
  Thus, $\eta=0$ and
$(\R^+\mu + \fg_\mu^\circ) \cap S=\R^+ \mu$.
\end{proof}

\begin{remark}  Proposition \ref{slice_reduction} allows us to always
 assume that, given that the hypothesis are satisfied, 
the element at which we are reducing is 
always fixed by the co-adjoint action of the symmetry group.  The hypothesis 
of Proposition \ref{slice_reduction} are satisifed for the action of any 
compact group. 
Indeed, if $G$ is compact, we can choose a $G-$invariant metric 
and hence an equivariant splitting $\fg^*=\fg_\mu^\circ \oplus \fg_\mu^*$, 
where $\fg_\mu^*$ is embedded as the normal fiber to $G \cdot \mu$ 
at $\mu$.  Take $S=\R^+ D_\epsilon$ where $D_\epsilon$ is a small 
$\epsilon$-ball about $0$ in $\fg_\mu^*$.  Then $S$ is a convex, 
$R^+$-invariant slice at $\mu$.
\end{remark}

We can now prove our second main theorem.

\begin{theorema}\label{reduction_theorem}Suppose a Lie group 
$G$ acts on a contact manifold 
$(M,\alpha)$, preseving $\alpha$ and let $\Phi_\alpha:M \to \fg^*$ be the 
associated
contact moment map.  Choose $\mu \in \fg^*$ and let $K_\mu$ be the connected
 Lie subgroup of 
$G_\mu$ with Lie algebra $\fk_\mu=\text{ker }\mu_{|\fg_\mu}$. Suppose
\begin{itemize} 
\item $K_\mu$ acts properly on $\Phi_\alpha\inv(\R^+\mu)$
\item A convex slice exists for $\mu$ which is invariant under 
dilations by $\R^+$
\end{itemize}
then the partition of the contact quotient by $G$
at $\mu$ by orbit types,

$$M_\mu=\coprod_{(H)} \frac{\Phi_\alpha\inv(\R^+ \mu) \cap M_{(H)}}{K_\mu}$$ is a 
stratification.  Here $M_{(H)}$ is the set of points whose 
stabilizer is conjugate to $H$ and the 
indexing set is the set of conjugacy classes of stabilizer subgroups of $K_\mu$.

\end{theorema}

\begin{proof}
  By Proposition \ref{slice_reduction}, we may assume that $\mu$ is fixed by
the co-adjoint action of $G$ on $\fg^*$.  Let $\Psi_\alpha:M \to \fk_\mu^*$ be the 
moment map for the action of $K_\mu$.  Because 
$\Psi_\alpha=i^T \circ \Phi_\alpha$,
where $i^T:\fg^* \to \fk_\mu^*$ is the natural projection, we have 
$\Psi_\alpha\inv(0)=\Phi_\alpha\inv(\R \mu)$.  Therefore 
$M_\mu$ is an open subset of $M/\!/K_\mu$, which is stratified by 
Theorem \ref{tscsq}. 
Open subsets of stratified spaces are naturally stratified.
\end{proof}

\section{Contact Structures on the Strata of the Reduced Space}
The previous section established the topological structure of the contact 
quotient  but did not address the geometrical structure of the quotient
under stratification.  
The point of this section is to show that the contact quotient is 
a contact stratified space; i.e, there exists a line bundle over $M_\mu$ 
which, when restricted to each stratum, defines a co-oriented 
contact structure.

\begin{proposition} Suppose a Lie group $G$ acts properly on a contact manifold 
$(M,\alpha)$, preserving $\alpha$, and let $\Phi_\alpha:M \to \fg^*$ be the 
associated
moment map. Let $H$ be a
isotropy subgroup of $G$, 
$$N=N(H)=\{g \in G|gHg\inv=H\}$$ the normalizer of $H$ in $G$, $L=N(H)/H$, and 
$$M_H=\{m \in M|G_m=H\}$$ the set of points of $M$ whose isotropy 
group is $H$.
Then $M_H$ is a contact submanifold of $M$, $L$ acts freely on $M_H$, and 
there is a diffeomorphism,
$$\vartheta:M_H/\!/L \to (\Phi_\alpha\inv(0) \cap M_{(H)})/G$$
\end{proposition}

\begin{proof} 
$M_H$ is a submanifold of $M$ and for all $x \in M_H$, 
$$T_x( M_H)=(T_xM)^H$$ where $(T_xM)^H$ is the set of $H$-fixed 
vectors in $T_xM$
(see, for example, Proposition 
27.5 of \cite{GS2}).
The Reeb vector field, $Y$, of $(M,\alpha)$, is 
$G$-invariant and, since the $H$-action preserves the contact
structure, 
$$(T_xM)^H=(\text{ker }\alpha_z \oplus \R Y(x))^H=
(\text{ker }\alpha_x)^H \oplus \R Y(x)$$
 Because $(\text{ker }\alpha_x)^H$ is a symplectic 
subspace of $(\text{ker }\alpha_x,d \alpha_{x|\text{ker }\alpha_x})$,
the restriction of $\alpha$ to $M_H$ is a contact 
form on $M_H$.

It what follows it is useful to cite Lemma 17, pg 220, of \cite{BL},
which identifies $\fl^*$, the dual of the Lie algebra of $L$, 
 with 
$(\fh^\circ)^H$.  
For any $X \in \fh$ and $x \in M_H$ $$\langle \Phi_\alpha(x), X \rangle 
=\alpha_x(X_M(x))=0$$ since $X_M(x)=0$.  
Hence, the image of $M_H$ in $\fg^*$ under 
$\Phi_\alpha$ is contained in  $\fh^\circ$.  Because the moment map is equivariant, 
$\Phi_\alpha(M_H) \subseteq (\fh^\circ)^H$.  
The action of $L$ on $M_H$ is defined by 
$nH \cdot x=n \cdot x$, where $nH$ is the coset containing $n \in N$.  This 
is free by definition.
The moment map, 
$\Psi_\alpha:M_H \to \fl^*$, for the $L$ action on $M_H$ is given by 
the restriction of $\Phi_\alpha$ to $M_H$.  Therefore, 
$$M_H/\!/L= (\Phi_\alpha\inv(0) \cap M_H)/L$$

The natural inclusion, 
$$\Phi_\alpha\inv(0) \cap M_H \to \Phi_\alpha\inv(0) \cap M_{(H)}$$
 descends to 
a map, $$\vartheta:M_H/\!/L \to (\Phi_\alpha\inv(0) \cap M_{(H)})/G$$ 
defined by $[x]_L \mapsto [x]_G$, where 
$[x]_L$ and $[x]_G$ denote the orbit through $x$ under the $L$ and $G$ actions 
respectively.  If $x \in \Phi_\alpha\inv(0) \cap M_{(H)}$, then the 
stablilizer of $x$ in $G$ is 
conjugate to $H$.  This implies that some element of the 
$G$-orbit through $x$ has stabilizer equal to $H$,
whence $\vartheta$ is surjective.

To show that $\vartheta$ is injective, suppose that $x,y \in \Phi_\alpha\inv(0) 
\cap M_H$ and that $y=g \cdot x$.  
Because $x$ and $y$ have stabilizer equal to $H$, 
it follows that $g \in N$ and therefore that $\vartheta$ is injective.

\end{proof}

By expressing each stratum of the contact quotient as the 
reduction at zero of a 
contact manifold by a freely acting 
symmetry group, we obtain a contact structure on 
each stratum.  It is slightly less clear, however, how these structures are 
related to one another.

\begin{remark} Suppose a Lie group $G$ acts on 
a manifold $N$ and $\sigma$ is an 
invariant 1-form on $N$.  
Then $\sigma$ is equivariant as a section, $N \to T^*N$, 
and hence descends to a section, 
$\bar \sigma:N/G \to (T^*N)/G$, where $G$ acts on 
$T^*N$ via the lifted action. 
\end{remark}

\begin{lemma}\label{stupid}  Suppose a Lie group $G$ acts on a line bundle 
$L \to X$ and $\sigma:X \to L$ is an invariant, non-vanishing 
section.  Then $L$ is equivariantly trivial and hence 
$L/G \to X/G$ is a vector bundle.
\end{lemma}

\begin{proof}  The trivialization $X \times \R \to L$ defined by 
$(x,t) \mapsto t \cdot \sigma(x)$ is equivariant (where $G$ 
acts trivially on $\R$) since $\sigma$ is.
\end{proof}

\begin{proposition} Suppose a 
 Lie group $G$ acts on a contact manifold 
$(M,\alpha)$, preserving $\alpha$, and let $\Phi_\alpha:M \to \fg^*$ be the 
associated 
moment map. Assume that $G$ acts locally freely and properly
on the zero level set of 
$\Phi_\alpha$.  Let $\xi=\text{ker }\alpha$ be the contact distribution 
on $M$ and $\xi^\circ$ its annihilator in $T^*M$.
Let $\xi^\circ/\!/G=(\xi^\circ_{|\Phi_\alpha\inv(0)})/G$.
Then there exists an embedding $\epsilon:\xi^\circ/\!/G \to T^*(M/\!/G)$ 
such that $\epsilon \circ \bar \alpha=\alpha_0$, where 
$\bar \alpha$ is the induced section and $\alpha_0$ is 
the reduced contact form 
on $M/\!/G$.
\end{proposition}

\begin{proof}  Set $Z=\Phi_\alpha\inv(0)$.
  Then $Z$ is a submanifold of $M$ and 
the natural 
inclusion, $i:TZ \to TM_{|Z}$ gives rise to a projection, 
$i^T:T^*M_{|Z} \to T^*Z$.   By Proposition \ref{basics}, 
the Reeb vector field is tangent to $Z$.
Split the tangent bundle of $M$ at $z$ 
as 
$$T_zM=\xi_z \oplus \R Y(z)$$ We obtain an induced splitting, 
$$T_zZ=(\xi_z  \cap T_zZ) \oplus \R Y(z)$$
If $0 \neq \eta \in \xi_z^\circ$, then 
$\langle i^T(\eta),Y(z) \rangle \neq 0$. Therefore, $i^T$ embeds 
$\xi^\circ$ into 
$T^*Z$.  Moreover, this embedding is equivariant.

Denote the symplectic moment map for the lifted $G$ action on 
$T^*Z$ by $\Psi$. Then 
$$\Psi\inv(0)=\{(z, \eta) \in T^*Z|\langle \eta,X_Z(z) \rangle=0 
\text{ for all } X \in \fg\}$$
For any $X \in \fg, z \in Z$, and $\eta \in \xi_z^\circ$, we have 
$$\langle i^T(\eta),X_Z(z) \rangle = \langle \eta,X_M(z) \rangle
=0$$ since $\Phi_\alpha(z)=0$ implies $X_M(z) \in \xi_z$ for 
all $X \in \fg$.  Hence, $\iota^T$ embeds $\xi^\circ_{|Z})$ into $\Psi\inv(0)$.
Let $\bar \iota^T:\xi^\circ/\!/G \to T^*(Z)/\!/G$ be the induced embedding.
The co-tangent bundle reduction theorem of Abraham-Marsden and Kummer 
(\cite{ma,ku}) asserts that there is a symplectomorphsim 
$\epsilon':T^*(Z)/\!/G \to T^*(M/\!/G)$.
Set $\epsilon=\epsilon' \circ \bar \iota^T$.

Recall that the reduced contact form on $M/\!/G$ is 
defined to be the unique 1-form, $
\alpha_0$, on $M/\!/G$ such that $\pi^*(\alpha_0)=\alpha_{|Z}$, 
where $\pi:Z \to 
M/\!/G$ is the orbit map.
It follows, by definition,
that $\pi^*(\epsilon \circ \bar \alpha)=\alpha_{|Z}$, giving 
$\epsilon \circ \bar \alpha=\alpha_0$.
\end{proof}

Since each stratum in the contact quotient at zero
can be expressed as the reduction of a contact manifold 
under a free group action, our third theorem follows immediately.

\begin{theorema}\label{coherence}
Suppose a Lie group $G$ acts properly on a 
contact manifold $(M,\alpha)$, preserving $\alpha$, and let $\Phi_\alpha:M \to 
\fg^*$ be the associated contact moment map.
For each stabilizer subgroup, $H$, of $G$, let 
$$(M/\!/G)_{(H)}:=(\Phi_\alpha\inv(0) \cap M_{(H)})/G$$ denote the stratum
associated
to $H$.  Then there exists a line bundle, $\xi^\circ/\!/G$ over 
$M/\!/G=\Phi_\alpha\inv(0)/G$ equipped 
with a section $\sigma$ such that 
for each such $H$, the restriction of $\xi^\circ/\!/G$
 and $\sigma$ to $(M/\!/G)_{(H)}$
defines a co-oriented contact  structure on the stratum. 
\end{theorema}

\begin{proof}  Set $Z:=\Phi_\alpha\inv(0)$ and $$\xi^\circ/\!/G:=(\xi^\circ_{|Z})/G$$
By Lemma \ref{stupid}, $\xi^\circ/\!/G$ is a line bundle over $M/\!/G$.
  Let $N$ be the normalizer of $H$ in $G$ and $L=N/H$.  Set 
$$M_H=\{m \in | G_m=H\}$$
 Then $M_H$ is a contact submanifold of $M$ and $L$ acts freely on
$M_H$, preserving $\alpha_{|M_H}$.  Denote the corresponding moment map by 
$\Psi_\alpha:M_H \to \fl^*$.  
Recall that we can identify $\Psi_\alpha\inv(0)$ with $\Phi_\alpha\inv(0) \cap M_H$
 and that there is a diffeomorphism 
$\vartheta:M_H/\!/L \to (M/\!/G)_{(H)}$.  It follows that, without being 
too fussy about equivalences versus equalities,  
$$(\xi^\circ/\!/G)_{|(M/\!/G)_{(H)}}=(\xi^\circ_{|\Phi_\alpha\inv(0) \cap
M_{(H)}})/G=(\xi^\circ_{|\Psi_\alpha\inv(0)})/L$$
Now apply the free case to obtain the result.
\end{proof}

Proposition \ref{slice_reduction} and Theorems \ref{reduction_theorem} and 
\ref{coherence} immediately imply 

\begin{corollary} \label{coherence_II}  Suppose a Lie group $G$ acts on 
a contact manifold $(M,\alpha)$, preserving $\alpha$, and let 
$\Phi_\alpha:M \to \fg^*$ be the associated moment map.  Choose $\mu 
\in \fg^*$ and suppose that 
\begin{itemize}
\item $K_\mu$ acts properly on $\Phi_\alpha\inv(\R^+\mu)$
\item A slice exists for $\mu$ which is invariant under dilations by $\R^+$
\item $\text{ker }\mu+\fg_\mu=\fg$
\end{itemize}
Let $K_\mu$ be the connected Lie subgroup of $G_\mu$ with Lie algebra 
$\fk_\mu=\text{ker }\mu_{|\fg_\mu}$.
For each stabilizer subgroup, $H$, of $G$, let 
$$(M_\mu)_{(H)}:=(\Phi_\alpha\inv(\R^+ \mu) \cap M_{(H)})/K_\mu$$ 
denote the stratum
associated
to $H$.  Then there exists a line bundle, $\xi^\circ_\mu$ over 
$$M_\mu=\Phi_\alpha\inv(\R^+\mu)/K_\mu$$ equipped 
with a section $\sigma$ such that 
for each such $H$, the restriction of $\xi^\circ_\mu$
 and $\sigma$ to $(M_\mu)_{(H)}$
defines a co-oriented contact  structure on the stratum.
\end{corollary}

\section{The Guillemin-Sternberg Quotient}

Let $G$ be a compact,
connected Lie group and $\Theta$ a submanifold of $\fg^*$ which is 
invariant under the co-adjoint action of $G$ on $\fg^*$.  Choose 
$\eta \in \Theta$ and let $\fh_\eta$ be the co-normal space of $\Theta$ at 
$\eta$.  Then $\fh_\eta$ is a Lie ideal in $\fg_\eta$.
 Let $H_\eta$ be the unique connected Lie subgroup of 
$G_\eta$ with Lie algebra $\fh_\eta$.  Call $\Theta$ {\bf proper} if 
$H_\eta$ is closed in $G_\eta$ for all $\eta \in \Theta$.

Suppose that $G$ acts in a Hamiltonian fashion on a 
symplectic manifold, $(N,\omega)$ and choose a corresponding equivariant 
moment map, $\Psi:N \to \fg^*$.  Let $\Theta$ be a proper, 
invariant submanifold of $\fg^*$.
If $\Psi$ is transverse to $\Theta$, then $Z:=\Psi\inv(\Theta)$ is a 
co-isotropic submanifold of $M$ and the 
leaf of the null foliation through $x \in Z$ is identified  
with the $H_{\Psi(x)}$ orbit through $x$.  

\begin{proposition} \label{gs_reduction}
Suppose a compact, connected Lie group $G$ acts in a 
Hamiltonian fashion on a symplectic manifold $(N,\omega)$ and let 
$\Psi:N \to \fg^*$ be a corresponding equivariant moment map.  Let
$\Theta \subset \fg^*$ is an invariant, proper submanifold of $\fg^*$ and 
set $Z:=\Psi\inv(\Theta)$.  If $\Psi$ is transverse to $\Theta$ then
the leaf space,
$$N_\Theta:=Z/\sim \hskip.25in \text{where } p \sim q \text{ if and only if }
p \text{ and } q \text{ are on the same leaf}$$
is a symplectic orbifold.
\end{proposition}

The  reader is refered to \cite{GS} for a proof of this Proposition.

\begin{example}  
If $\Theta=0 \in \fg^*$, then $N_\Theta$ is the 
Marsden-Weinstein-Meyer reduced space.  If $\Theta=\O$, a co-adjoint orbit, 
then $N_\Theta$ is the Kahzdan-Kostant-Sternberg reduced space.
\end{example}

\begin{definition}
A symplectic manifold  $(N,\omega)$ is called a {\bf symplectic cone} 
if there is a free,
proper $\R$ action $\tau:N \times \R \to N$ such that $\tau_t^* \omega=e^t \omega$.
\end{definition}

\begin{example}\label{can_action}
  Let $(M,\alpha)$ be a compact, co-oriented contact manifold 
and $(M \times \R,d(e^t \alpha))$ the symplectization of $(M,\alpha)$.  Then the
action of $\R$ on $M \times \R$ given by $s \cdot (m,t)=(m,t+s)$ makes 
$(M \times \R,d(e^t \alpha))$ into a symplectic cone.
\end{example}

In fact, one can show that all symplectic cones are of the form given in the above 
example.

\begin{proposition}\label{can_form}  Let $(N,\omega)$ be a symplectic cone.  Then 
$N \simeq M \times \R$ where $M=N/\R$ is a co-oriented contact manifold.
\end{proposition}

The reader is refered to, for example, \cite{Le} for a proof of the above 
proposition and more 
details on symplectic cones.

\begin{definition}  Call $\mu \in \fg^*$ 
{\bf integral} if there is a group homomorphism $\chi_\mu:G_\mu \to S^1$ such that
$d \chi_\mu=\mu_{|\fg_\mu}$.  Call a co-adjoint 
orbit $\O$ an {\bf integral} orbit if it is an orbit through an integral 
element of $\fg^*$.
\end{definition}

\begin{remark}
Note that if $\mu$ is integral, then $H_\mu=\text{ker }\chi_\mu=K_\mu$, the 
kernel group of $\mu$ defined earlier.  Integral orbits are automatically 
proper.  The cone on an integral orbit $\R^+ \O$ is a proper submanifold of 
$\fg^*$.
\end{remark}

Let $(M,\alpha)$ be a co-oriented contact manifold equipped with an action of a
compact, connected Lie group $G$ which preserves $\alpha$ and let 
$\Phi_\alpha:M \to \fg^*$ be the associated contact moment map.
Extend the $G$ action to the symplectization $(M \times \R,d(e^t \alpha))$ by 
$g \cdot (m,t)=(g \cdot m,t)$.  A corresponding
equivariant symplectic moment map $$\Psi:M \times \R \to \fg^*$$
 is given by $\Psi(m,t)=e^t \Phi_\alpha(m)$.  
Let $\O \subset \fg^*$ be an integral 
orbit and assume that $\Psi$ is transveral to $\R^+ \O$.  
By Proposition \ref{gs_reduction},
$(M \times \R)_{\R^+ \O}$ is a symplectic orbifold.  The $\R$ action on 
$M \times \R$ described 
in Example \ref{can_action} descends to $(M \times \R)_{\R^+ \O}$.
By Proposition \ref{can_form},   the orbit space under the 
$\R$ action on $(M \times \R)_{\R^+ \O}$ is thus a contact manifold and 
is called the 
{\bf Guillemin-Sternberg} reduction of $M$ by $\O$. It is denoted by 
$M^{GS}_\O$.

The Guillemin-Sternberg reduction should be thought of as a larger form of 
contact reduction in the following sense. Since $\R^+ \O$ is $\R^+$ invariant, 
$\Psi\inv(\R^+ \O)=\Phi_\alpha\inv(\R^+ \O) \times \R$.  For any $\mu \in \O$
the contact quotient $M_\mu$ is a contact suborbifold of $M^{GS}_\O$.
More, in fact, is true.
Since the leaf of the null foliation through $(p,t)$ is given by the 
$K_{\Phi_\alpha(p)}$ orbit through $(p,t)$, the restriction of $\Psi$ to 
a leaf in $\Psi\inv(\R^+ \O)$ is constant.  Therefore, $\Psi$ descends to
$$\tilde \Psi:(M \times \R)_{\R^+ \O} \to \R^+ \O$$
For any $t \in \R$, let $\epsilon^\mu_t:\Phi_\alpha\inv(\R^+ \mu) \to \R$ be the map 
$m \mapsto t-t_m$, where $\Phi_\alpha(m)=e^{t_m} \mu$.  Denote the graph of 
$\epsilon_t^\mu$ by $\Gamma_t^\mu$.  It follows that $\Gamma_t^\mu=\Psi\inv(e^t \mu)$ 
and, since $K_{e^s \mu}=K_{e^t \mu}$, that 
$$\tilde \Psi\inv(e^t \mu)=\Gamma_t^\mu/K_\mu \simeq M_\mu$$
In summary, we have shown,

\begin{proposition}  Let $(M,\alpha)$ be a compact, connected, co-oriented
contact manifold equipped with an action of a compact, connected Lie group $G$ 
which preserves $\alpha$ and let $\Phi_\alpha:M \to \fg^*$ be the associated 
moment map.  Choose an integral $\mu \in \fg^*$ and let $\O$ be the 
co-adjoint orbit through $\mu$.  Then $M_\O^{GS}$ fibers over $\O$ with 
typical fiber $M_\mu$.
\end{proposition}

\appendix
\section{Albert's Quotient}
Independent of Guillemin-Sternberg, C. Albert discovered 
an elegant method for contact reduction.  However, his method 
depends upon the choice of contact form used to represent the given 
contact structure.  Let $(M,\alpha)$ be a co-oriented contact 
manifold on which a compact group $G$ acts by contact transformations.
Let $\Phi_\alpha:M \to \fg^*$ be the associated moment map.
Suppose $\mu$ is a regular value of $\Phi_\alpha$, so that
$Z:=\Phi_\alpha\inv(\mu)$  is a submanifold of 
$M$.  Consider the map, $\tilde \tau: \fg_\mu \to \chi(Z)$ defined by 
$A \mapsto A_Z-\langle \mu,A \rangle Y$, where $Y$ is the 
Reeb vector field.
It follows that $\tilde \tau$ is a map of Lie algebras.  
By a 
theorem of Palais
 there is a unique map, $\tau:\tilde G_\mu \to \text{Diff}(Z)$, 
where $\tilde G_\mu$ is the universal covering group of $G_\mu$,
whose differential is $\tilde \tau$.  Set $H=\tilde G_\mu/\text{ker }\tau$.  
Then 
$H$ acts locally freely and effectively 
on $Z$.  The {\bf Albert reduction} of $M$ by $G$ at 
$\mu$ is defined to be $$M^A_\mu:=Z/H$$  If the $H$ action is proper, 
this is an orbifold and, 
by definition, $\alpha$ descends to $\alpha_\mu$ on $M^A_\mu$.
A standard argument shows that $\alpha_\mu$ is a contact form.
However $f \alpha$, where $f$ is a positive 
invariant function, is an invariant contact form with associated 
moment map $\Phi_{f \alpha}=f \Phi_\alpha$.  This indicates that 
the Albert quotient will be dependent on the choice of contact 
form.  The following example illustrates this
dependence.

\begin{example}
Let 
$$\aligned
\E_1&=\{(z_1,z_2,z_3) \in \C^3||z_1|^2+2 |z_2|^2+2|z_3|^2=1\} \cr
\E_2&=\{(z_1,z_2,z_3) \in \C^3|\frac{1}{2}|z_1|^2+ |z_2|^2+
\frac{1}{3}|z_3|^2=1\}
\endaligned$$

Let $S^1$ act on $\E_i$, $i=1,2$, by the restriction of the action on 
$\C^3$ given 
by the weights $(1,-1,-1)$.  Since both $\E_i$ are star shaped about 
the origin in $\C^3$, 
they are isomorphic contact manifolds.  The difference in shape amounts 
to a choice 
of different contact forms in the same conformal class.
Let $\Phi_i$ denote the contact moment map for the $S^1$ action on $\E_i$.
A simple calculation yields
$$ 
\Phi_1(z_1,z_2,z_3)=\Phi_2(z_1,z_2,z_3)=|z_1|^2-|z_2|^2-|z_3|^2$$

It is not hard to see that 
$$\Phi_1\inv(1)=\{(z_1,0,0) \in \C^3||z_1|^2=1\}$$
and that the Albert action is trivial.  Hence the reduction of $\E_1$ is 
the circle, $S^1$.

A more tedious calculation gives

$$\aligned 
\Phi_2\inv(1)=\{(z_1,z_2,z_3) \in \C^3|&|z_1|^2=\frac{2}{3}
(2+\frac{2 |z_3|^2}{3}), 
|z_2|^2=-\frac{2}{3}(-\frac{1}{2}+\frac{5 |z_3|^2}{6}),\cr
&0 \leq |z_3|^2 \leq \frac{3}{5}\} \endaligned$$

The symmetry group in the Albert quotient is the circle, acting with weights 
$(\frac{1}{2},-2,-\frac{4}{3})$.
 Note that 
$\Phi_2\inv(1)$ contains a 3-torus, namely 
$$\{(z_1,z_2,z_3) \in \C^3||z_1^2=\frac{14}{9},|z_2|^2=
\frac{1}{18},|z_3|^2=\frac{1}{2}\}$$
Hence, the Albert reduction of $\E_2$ is 
at least 3 dimensional, illustrating the depending of the Albert 
quotient on the 
choice of the contact form.
\end{example}

\end{document}